# Chance-Constrained Two-Stage Unit Commitment under Uncertain Load and Wind Power Output Using Bilinear Benders Decomposition

Yao Zhang, *Student Member, IEEE*, Jianxue Wang, *Member, IEEE*, Bo Zeng, *Member, IEEE*, and Zechun Hu, *Member, IEEE*

*Abstract*—In this paper, we study unit commitment (UC) problems considering the uncertainty of load and wind power generation. UC problem is formulated as a chance-constrained two-stage stochastic programming problem where the chance constraint is used to restrict the probability of load imbalance. In addition to the conventional mixed integer linear programming formulation using Big-M, we present the bilinear mixed integer formulation of chance constraint, and then derive its linear counterpart using McCormick linearization method. Then, we develop a bilinear variant of Benders decomposition method, which is an easy-to-implement algorithm, to solve the resulting large-scale linear counterpart. Our results on typical IEEE systems demonstrate that (i) the bilinear mixed integer programming formulation is stronger than the conventional one; (ii) the proposed Benders decomposition algorithm is generally an order of magnitude faster than using a professional solver to directly compute both linear and bilinear chance-constrained UC models.

*Index Terms*—Benders decomposition, chance constraint, unit commitment, bilinear formulation, stochastic programming, wind power.

## Nomenclature

### A. Indices and Sets

$b \in \mathcal{B}$    Index of buses, from 1 to $B=|\mathcal{B}|$. The notation $|.|$ represents the cardinality of a set.

$g \in \mathcal{G}$    Index of thermal units, from 1 to $G=|\mathcal{G}|$.

$l \in \mathcal{L}$    Index of transmission lines, from 1 to $L=|\mathcal{L}|$.

$n \in \Omega$    Index of scenarios, from 1 to $N=|\Omega|$.

$q \in \mathcal{Q}$    Index of wind farms, from 1 to $Q=|\mathcal{Q}|$.

$s \in \mathcal{S}_g$    Index of startup segments, from 1 to $S_g=|\mathcal{S}_g|$.

$t \in \mathcal{T}$    Index of time periods, from 1 to $T=|\mathcal{T}|$.

$\mathbf{G}(b), \mathbf{W}(b)$    Thermal units/wind farms subset located at bus $b$.

### B. Variables

1) *Continuous and Non-negative Variables:*

$P_{gt}$    Power output of thermal unit $g$ at period $t$ [MW].

$P_{gtn}$    Power output of thermal unit $g$ at period $t$ under scenario $n$ [MW].

$P_{gt}(\xi)$    Power output of thermal unit $g$ at period $t$ for the second stage [MW].

$R_{gt}^{\pm}$    Up/down spinning reserve provided by thermal unit $g$ at period $t$ [MW].

$r_{gtn}^{\pm}$    Up/down spinning reserve deployed by thermal unit $g$ at period $t$ under scenario $n$ [MW].

$r_{gt}^{\pm}(\xi)$    Up/down spinning reserve deployed by thermal unit $g$ at period $t$ for the second stage [MW].

$\eta_{tn}^{1\pm}$    Slack variable for the power balance constraint at period $t$ under scenario $n$.

$\eta_{ltn}^{2}$    Slack variable for the capacity constraint of transmission line $l$ at period $t$ under scenario $n$.

2) *Binary Variables:*

$u_{gt}$    Commitment status that is equal to 1 if thermal unit $g$ is online at period $t$.

$v_{gt}$    Startup status that is equal to 1 if thermal unit $g$ starts up at period $t$.

$y_{gt}$    Shutdown status that is equal to 1 if thermal unit $g$ shuts down at period $t$.

$z_n$    Binary indicator that is equal to 1 if the associated constraint is not satisfied under scenario $n$.

$\delta_{gst}$    Startup type $s$ of thermal unit $g$, which is equal to 1 when the unit starts up at the period $t$ and has been offline within $[T_{gs}^{SU}, T_{g,s+1}^{SU})$ hours.

### C. Parameters

$C_g^{NL}$    No-load cost of thermal unit $g$ [\$].

$C_g^{R\pm}$    Up/down spinning reserve cost of thermal unit $g$ [\$].

$C_g^{SD}$    Shutdown cost of thermal unit $g$ [\$].

$C_{gs}^{SU}$    Startup cost of thermal unit $g$ under startup type $s$ [\$].

$C_l$    Capacity of transmission line $l$ [MW].

$K_{lb}$    Power flow distribution factor for the transmission line $l$ due to the net injection at bus $b$.

$L_{bt}$    Forecasted load demand located at bus $b$ at period $t$ [MW].

$L_{btn}$    Load demand located at bus $b$ at period $t$ under scenario $n$ [MW].

$L_{bt}(\xi)$    Random parameter indicating the uncertain load demand located at bus $b$ at period $t$ [MW].

$P_g^{max}, P_g^{min}$    Maximum/minimum power output of thermal unit $g$ [MW].

$R_g^{\pm max}$    Maximum amount of up/down spinning reserve capability provided by thermal unit $g$ [MW].

$RU_g, RD_g$    Ramp up/down rate of thermal unit $g$ [MW/h].

This work was supported in part by the National Natural Science Foundation of China (51277141). Y. Zhang and J. Wang are with the School of Electrical Engineering, Xi'an Jiaotong University, Xi'an, Shaanxi, 710049, China (e-mail: zy.3110161019@stu.xjtu.edu.cn; jxwang@mail.xjtu.edu.cn). B. Zeng is with the Department of Industrial Engineering, University of Pittsburgh, Pittsburgh, PA 15260 USA (e-mail: bzeng@pitt.edu). Zechun Hu is with the Department of Electrical Engineering, Tsinghua University, Beijing, 100084, P. R. China (e-mail: zechhu@tsinghua.edu.cn).

| | |
|---|---|
| $SU_g, SD_g$ | Startup/shutdown capability of thermal unit $g$ [MW]. |
| $TU_g, TD_g$ | Minimum uptime/downtime of thermal unit $g$ [h]. |
| $W_{qt}$ | Forecasted power output of wind farm $q$ at period $t$ [MW]. |
| $W_{qtn}$ | Power output of wind farm $q$ at period $t$ under scenario $n$ [MW]. |
| $W_{qt}(\xi)$ | Random parameter indicating the uncertain power output of wind farm $q$ at period $t$ [MW]. |
| $\pi_n$ | Probability of the scenario $n$. |

## I. INTRODUCTION

Due to low cost and low emission, renewable power generation, e.g., wind power generation, has developed rapidly all over the world in recent decades. Nevertheless, wind power generation is intermittent and it is quite difficult to give an accurate day-ahead prediction [1]. As a conventional method in power system operation, the uncertainty caused by load variation and generator's forced outage is handled by imposing predefined reserve requirements. This method, which is known as reserve adjustment method, is easy to implement in practice and has been widely adopted in today's power industry for many years [2]. However, dispatching extra generators as reserves is uneconomic to deal with the uncertainty, especially when the reserve requirement is determined by some rather simple rules [3]. Indeed, the volatility from wind power generation is generally higher than load variation. Thus, even with reserves determined by rules, power systems may still suffer from the reserve scarcity when wind power output deviates significantly from the predicted value [4]. So, to more economically and reliably improve wind power penetration level in the grid, many research efforts have been devoted to using advanced methods to integrate power system operation, especially unit commitment (UC) [5], with analytically described wind generation randomness. Among those methods, two most popular approaches are stochastic programming [3, 6-14] and robust optimization [4, 15-19].

The most general approach in stochastic programming is to utilize a set of representative scenarios (through sampling if necessary) to capture random factors, e.g., wind power generation [7], and introduce a recourse decision problem for every scenario. As a result, a deterministic UC model will be converted into a two-stage stochastic UC (SUC) model [9]. This model minimizes the expected cost while satisfying operational constraints under those scenarios [10]. Thus, it guarantees that commitment decisions of conventional generators are sufficiently flexible to address the uncertainty associated with wind power generation [12]. Compared with deterministic UC (DUC) models, SUC models have advantages of high reliability, as shown in [3] and [11]. Nevertheless, because all scenarios must be considered, which may include some extreme scenarios, computing SUC could lead to costly solutions.

In fact, enforcing a complete coverage over all possible extreme scenarios could be physically and economically impractical. To address that issue, chance-constrained optimization is introduced to restrict the consideration of rarely occurred extreme scenarios. Specifically, a small number of scenarios, whose realization probabilities sum up to $\epsilon$ (<1), can be ignored in deriving an optimal solution. Clearly, when $\epsilon$ equals 0, the chance constraint disappears and all scenarios must be considered. By adjusting the value of $\epsilon$, the decision maker will be able to have a desired trade-off between cost and reliability.

Because of the aforementioned advantages, the chance-constrained optimization has been employed as decision tools for UC problems in the last decade [20]. In [21-23], chance-constrained UC models have been developed where a large portion of wind power output should be utilized with high probability. In the chance-constrained UC problem reported by [24], reserve requirements and transmission line capacity limitations are formulated as chance constraints to maintain the reliability of system operation. In addition to UC problems, chance-constrained optimization method has also been employed in other research fields, such as reserve scheduling [25], generation expansion planning [26], transmission expansion planning [27], and demand response [28].

Although chance-constrained optimization is powerful in modeling, it is well recognized that solving chance-constrained optimization problems is computationally challenging, especially when the random variable $\xi$ follows an unstructured and continuous distribution. Note that in general, analytically representing the probability constraint $Pr\{G(x,\xi) \leq \mathbf{0}\}$ (where $x$ represents the set of decision variables and $\xi$ denotes a random realization) is very difficult as it requires the complex computation of multi-dimensional integration. Certainly, for some special cases, chance constraint can be converted into a closed form expression and the whole model can be reformulated into a regular mixed integer programming (MIP) model (e.g., [20, 24, 26, 29]). In [30, 31], such property is used to study chance constrained optimal power flow problems subject to uncertain parameters of the normal distribution.

Nevertheless, because general cases, especially those with joint chance constraints, are complicated, by following the convention in stochastic programming, sampling-based method has been extended to generate a finite set of scenarios and then chance constraints are imposed upon those sampled scenarios. Specifically, by using a binary variable to indicate whether the associated scenario should be considered (i.e., the corresponding constraints must be satisfied) or ignored (i.e., the corresponding constraints can be violated), the chance-constrained model can be converted into an MIP model (through using Big-M method). Such strategy has been adopted in [21-23, 25], where, however, professional MIP solvers are often incapable to compute the resulting large-scale MIPs and fast heuristic methods are necessary [22, 25]. Hence, it remains a critical challenge to efficiently compute chance-constrained UC problems for real applications, especially when a large number of scenarios are needed and joint chance constraints are required.

We note that a formulation developed based on the concept of *conditional value-at-risk* (CVaR) sometimes is also considered as a chance-constrained model [32]. As pointed in [5], although both chance-constrained and CVaR-based UC models reflect the modeler's risk consideration, the former one is computationally much more challenging. Indeed, chance constraint usually introduces non-convexity into the original model while CVaR preserves convexity [33].

In this paper, we aim to develop an efficient algorithm to address the computational challenge associated with chance-constrained UC problem. Specifically, a chance-constrained two-stage UC model is first presented where the power imbalance due to random load demand and wind power generation is restricted by a predefined low probability. Then, we adopt the bilinear reformulation to represent chance constraints [34], and convert that bilinear model into a linear one through McCormick linearization method [35]. To deal with the large number of stochastic scenarios, we provide the bilinear Benders reformulation of the original model and customize the bilinear variant of Benders decomposition algorithm to achieve fast computation [34]. Our numerical results on standard IEEE systems, including 118-bus system, demonstrate that (i) the bilinear UC formulation is stronger than the widely adopted Big-M based formulation; (ii) our proposed bilinear Benders decomposition algorithm is generally an order of magnitude faster than using a professional solver to directly compute both linear and bilinear chance-constrained UC models. To the best of our knowledge, this is the first effort to develop a fast algorithm to solve large-scale chance-constrained UC problems. Also, we believe that, a great computational improvement can be achieved in solving existing UC variants in [21-23] by adopting this bilinear Benders decomposition method.

The rest of this paper is organized as follows. A chance-constrained two-stage UC problem is formulated in Section II. Section III introduces the bilinear reformulation of chance constraints and the bilinear Benders decomposition algorithm. Section IV provides numerical results from case studies using an illustrative six-bus system and a modified 118-bus system. This paper is concluded in Section V.

## II. MATHEMATICAL FORMULATION

In this section, we present a unit commitment problem that is formulated as a two-stage model. It extends the stochastic UC model with non-recourse cost in [36] by imposing joint chance constraints on random wind power output and load demand. As our main purpose is to develop efficient algorithms to meet the computational needs of chance-constrained UC problems, e.g., [21-23], we do not consider the uncertainty of component outages. As mentioned earlier, more advanced algorithm development should be necessary to improve our solution capacity from CVaR based contingency-constrained UC models [32] to actual chance constraint based ones.

The objective function (1) aims to minimize the day-ahead operational cost. It consists of fuel cost, startup/shutdown cost and up/down spinning reserve cost of all thermal units over the entire scheduling periods. The first stage (constraints (2)-(14)) refers to the day-ahead scheduling in the commitment and the dispatch of thermal units considering short-term deterministic (point) prediction of wind power generation. The second stage (constraints (15)-(23)) refers to the redispatching of thermal units for satisfying all operational constraints under the uncertainty of wind power generation and load demand (represented by a random variable $\xi$). The first and second stages are known as preventive and corrective stages, respectively [36]. Also, this formulation considers up/down spinning reserve deployments and reserve deliverability across the transmission network [37]. The complete formulation of chance-constrained UC problem is described as follows.

- Objective function

$$\min \sum_{t \in \mathcal{T}} \sum_{g \in \mathcal{G}} [C_g^{\text{NL}} u_{gt} + \sum_{s \in \mathcal{S}_g} C_{gs}^{\text{SU}} \delta_{gst} + C_g^{\text{SD}} y_{gt} + F_g(P_{gt}) + C_g^{\text{R+}} R_{gt}^+ + C_g^{\text{R-}} R_{gt}^-] \quad (1)$$

- Constraints for the first stage

$$\text{s.t.} \quad \delta_{gst} \leq \sum_{i=T_{gs}^{SU}}^{T_{g,s+1}^{SU}-1} y_{g,t-i} \quad \forall g, t, s = 1, \dots |\mathcal{S}_g| - 1 \quad (2)$$

$$\sum_{s \in \mathcal{S}_g} \delta_{gst} = v_{gt} \quad \forall g, t \quad (3)$$

$$\sum_{i=t-TU_g+1}^{t} v_{gi} \leq u_{gt} \quad \forall g, t \quad (4)$$

$$\sum_{i=t-TD_g+1}^{t} y_{gi} \leq 1 - u_{gt} \quad \forall g, t \quad (5)$$

$$u_{gt} - u_{g,t-1} = v_{gt} - y_{gt} \quad \forall g, t \quad (6)$$

$$\sum_{g \in \mathcal{G}} P_{gt} + \sum_{q \in \mathcal{Q}} W_{qt} = \sum_{b \in \mathcal{B}} L_{bt} \quad \forall t \quad (7)$$

$$-C_l \leq \sum_{b \in \mathcal{B}} K_{lb} \left( \sum_{g \in \mathbf{G}(b)} P_{gt} + \sum_{q \in \mathbf{W}(b)} W_{qt} - L_{bt} \right) \leq C_l \ \forall l, t \quad (8)$$

$$P_{gt} + R_{gt}^+ - P_{g,t-1} \leq RU_g u_{g,t-1} + SU_g v_{gt} \quad \forall g, t \quad (9)$$

$$-P_{gt} + R_{gt}^- + P_{g,t-1} \leq RD_g u_{gt} + SD_g y_{gt} \quad \forall g, t \quad (10)$$

$$P_{gt} + R_{gt}^+ \leq P_g^{\max} u_{gt} \quad \forall g, t \quad (11)$$

$$P_{gt} - R_{gt}^- \geq P_g^{\min} u_{gt} \quad \forall g, t \quad (12)$$

$$0 \leq R_{gt}^+ \leq R_g^{+\max} \quad \forall g, t \quad (13)$$

$$0 \leq R_{gt}^- \leq R_g^{-\max} \quad \forall g, t \quad (14)$$

- Constraints for the second stage

$$\Pr\left( \sum_{g \in \mathcal{G}} P_{gt}(\xi) + \sum_{q \in \mathcal{Q}} W_{qt}(\xi) = \sum_{b \in \mathcal{B}} L_{bt}(\xi), \forall t \right) \geq 1 - \epsilon \quad (15)$$

$$-C_l \leq \sum_{b \in \mathcal{B}} K_{lb} \left[ \sum_{g \in \mathbf{G}(b)} P_{gt}(\xi) + \sum_{q \in \mathbf{W}(b)} W_{qt}(\xi) - L_{bt}(\xi) \right] \leq C_l$$

$$\forall l, t \quad (16)$$

$$P_{gt}(\xi) = P_{gt} + r_{gt}^+(\xi) - r_{gt}^-(\xi) \quad \forall g, t \quad (17)$$

$$0 \leq r_{gt}^+(\xi) \leq R_{gt}^+ \quad \forall g, t \quad (18)$$

$$0 \leq r_{gt}^-(\xi) \leq R_{gt}^- \quad \forall g, t \quad (19)$$

$$P_{gt}(\xi) - P_{g,t-1}(\xi) \leq RU_g u_{g,t-1} + SU_g v_{gt} \quad \forall g, t \quad (20)$$

$$-P_{gt}(\xi) + P_{g,t-1}(\xi) \leq RD_g u_{gt} + SD_g y_{gt} \quad \forall g, t \quad (21)$$

$$P_{gt}(\xi) \leq P_g^{\max} u_{gt} \quad \forall g, t \quad (22)$$

$$P_{gt}(\xi) \geq P_g^{\min} u_{gt} \quad \forall g, t \quad (23)$$

In the above formulation, $F_g(.)$ represents the quadratic fuel cost function which can be approximated by a piecewise linearization. The startup cost of thermal unit is approximated by a stair-wise function [38, 39]. Constraints (2) and (3) choose the suitable startup-type variable $\delta_{gst}$ that activates the corresponding startup cost $C_{gs}^{SU}$ in the objective function [40]. Constraints (4) and (5) describe the minimum uptime and downtime limitations of all thermal units. Constraint (6) guarantees that binary variables $v_{gt}$ and $y_{gt}$ get proper values at startup and shutdown times. Constraints (7) and (8) indicate power balance constraints and transmission line capacity constraints, respectively. Constraints (9) and (10) make sure that the thermal unit operates within the ramping rate limitation between two successive periods. Constraints (11) and (12) restrict the minimum and maximum output of thermal units. Constraints (13) and (14) restrict the amount of up/down spinning reserve capability provided by each thermal unit. The above constraints belong to the first stage formulation considering deterministic forecasts. In the second stage constraints, the uncertainty of wind power generation and load consumption is characterized by random variable $\xi$. The chance constraint (15) guarantees that the probability of load imbalance should be less than a predefined risk level $\epsilon$. It is utilized to control the power-demand balance and avoid the load curtailment. Besides, the second stage constraints also include transmission line capacity constraints (16), actual power output constraints of thermal units (17), up/down spinning reserve deployment constraints (18)-(19), unit ramping rate constraints (20)-(21) and unit generation capacity constraints (22)-(23) [37].

## III. SOLUTION METHODOLOGY

In this section, to compute the aforementioned chance-constrained UC model, we present a non-traditional bilinear MIP reformulation of chance constraints. Then, we develop a bilinear variant of Benders decomposition algorithm. Following the convention of sampling based methods, the general and continuous distribution of random variable is represented by a set of finite and discrete scenarios through the Monte Carlo sampling method, which converts the chance-constrained UC model into a computationally friendly form. Next, we describe three parts of our whole scheme: A) scenario generation; B) bilinear reformulation, and C) bilinear Benders decomposition.

### A. Scenario Generation

In this paper, Monte Carlo simulation is used to generate scenarios for wind power output and load consumption. It usually requires assuming probability distribution of the uncertainty. However, this assumption may be unrealistic because it is difficult to accurately identify the shape of uncertainty distribution for day-ahead UC problems. In this paper, probabilistic forecasting is developed to predict probabilistic information (e.g., probabilistic distribution) of the hourly wind farm production and load consumption for the next 24 or 48 hours. The used probabilistic forecasting approach is non-parametric without any assumption of the density shape [41, 42]. In comparison with the conventional approaches, the novel approach makes UC solutions immune against the parameter uncertainty of probabilistic distribution. In addition, its predicted information dynamically varies with the future weather condition (provided by numerical weather prediction). Hence, such probabilistic forecasting provides more accurate probabilistic information of the uncertainty than the traditional approaches.

A general framework of probabilistic forecasting for renewable energy generation, especially for wind power generation, has been presented in previous articles [43] and [44]. It is a combination of *k*-nearest neighbors (*k*-NN) algorithm and kernel density estimator (KDE) method. The effectiveness of this approach has been testified through Global Energy Forecasting Competition 2014 (GEFCom 2014, http://www.gefcom.org). In addition, probabilistic load forecasting proposed in [45] is utilized to produce the information of probabilistic distribution for the future load consumption. Then, $N$ scenarios including wind power output and load demand are generated from probabilistic forecasting [46]. Each scenario has the same probability $\pi_n = 1/N$. Finally, the second stage random variables $W_{qt}(\xi)$ and $L_{bt}(\xi)$ are replaced by their scenarios $W_{qtn}$ and $L_{btn}, 1 \leq n \leq N$, respectively. Other second stage decision variables are represented by $P_{gtn}, r_{gtn}^+$ and $r_{gtn}^-, 1 \leq n \leq N$.

### B. Bilinear Reformulation of Chance Constraint

Traditionally, by introducing a binary indicator variable and using Big-M method, the chance constraint (15) can be equivalently reformulated with a set of linear constraints as follows [34]:

$$-Mz_n \leq \sum_{g \in \mathcal{G}} P_{gtn} + \sum_{q \in \mathcal{Q}} W_{qtn} - \sum_{b \in \mathcal{B}} L_{btn} \leq Mz_n \quad \forall t, n \quad (24)$$

$$\sum_{n=1}^{N} \pi_n z_n \leq \epsilon, \ z_n \in \{0,1\}, \ \forall n \quad (25)$$

where $M$ is a sufficiently large number. Note that $z_n$ is a binary indicator: $z_n = 0$ indicates a responsive scenario where the power balance constraints in (24) must be met for all $t$, and $z_n = 1$ indicates a non-responsive scenario where the constraints in (24) could be ignored (i.e., load imbalance), due to the Big-M parameter. Then, the constraint (25) is introduced as an equivalence of chance constraint to restrict the number of $z_n, 1 \leq n \leq N$, being ones (i.e., restrict the number of non-responsive scenarios). As a result, the chance-constrained UC problem is converted into an MIP problem that can be readily computed by any MIP solvers. This solution strategy is adopted in computing chance-constrained UC models within the sample averaging approximation (SAA) framework [21-23].

Nevertheless, such Big-M formulation has serious issues: (i) it is challenging to estimate a reasonable value for Big-M without hurting the equivalence between constraints (15) and (24); (ii) Big-M always slows down the computational efficiency, rendering the large-scale instances practically unsolved. Different from that conventional Big-M MIP reformulation, we adopt the strategy developed in [34] to derive another equivalent reformulation, i.e., a bilinear mixed integer reformulation of chance constraint as the following:

$$(\sum_{g \in \mathcal{G}} P_{gtn} + \sum_{q \in Q} W_{qtn} - \sum_{b \in \mathcal{B}} L_{btn})(1 - z_n) = 0 \quad \forall t, n \quad (26)$$

$$\sum_{n=1}^{N} \pi_n z_n \leq \epsilon, \ z_n \in \{0,1\}, \ \forall n \quad (27)$$

The effectiveness of bilinear reformulation (26) is obvious: the non-responsive scenario will be removed from the whole formulation by assigning $z_n$ to be one. Although the proposed bilinear constraint (26) is nonlinear, it can be easily converted into linear constraints by McCormick linearization method [35] as the following. The bilinear term $P_{gtn} z_n (= \tilde{P}_{gtn})$ is linearized through including the constraints (29)-(31).

$$\sum_{g \in \mathcal{G}} (P_{gtn} - \tilde{P}_{gtn}) + (\sum_{q \in Q} W_{qtn} - \sum_{b \in \mathcal{B}} L_{btn})(1 - z_n) = 0 \forall t, n (28)$$

$$0 \leq \tilde{P}_{gtn} \leq P_g^{\max} z_n \quad (29)$$

$$\tilde{P}_{gtn} \leq P_{gtn} \quad (30)$$

$$\tilde{P}_{gtn} \geq P_{gtn} - P_g^{\max}(1 - z_n) \quad (31)$$

**Remark:** After McCormick linearization, we obtain a linear MIP reformulation. As shown in Section 2.1 of [47], because $z_n$ is binary and $P_g^{\max}$ is the upper bound of $P_{gtn}$, the formulation defined by (29)-(31) is ideal, i.e., its linear programming relaxation is the convex hull of $\{\omega = P_{gtn} z_n, 0 \leq P_{gtn} \leq P_g^{\max}, z_n \in \{0,1\}\}$. Thus, no other formulation strongly dominates this one. Given that $P_{gtn}$ is naturally and tightly bounded by the maximum power output $P_g^{\max}$ and Big-M parameter is typically set to be large, this McCormick linearized reformulation should be theoretically stronger than the conventional Big-M linear formulation.

### C. Bilinear Benders Decomposition Algorithm

Although the chance-constrained two-stage UC problem can be reformulated as an MIP, this problem becomes large as the increase of scenarios, which cannot be solved by off-the-shelf MIP solvers with efficiency. In this subsection, on top of our bilinear reformulation, a bilinear variant of Benders decomposition [34] is designed to efficiently solve chance-constrained UC problems. Because the objective function (1) only contains the day-ahead operational cost, we only need to check the feasibility of the second-stage subproblem under each scenario. It needs no return of optimality cuts from subproblems to the master problem. As a result, the presented Benders decomposition algorithm is employed to decompose this problem into master UC problem and several feasibility check subproblems corresponding to all scenarios. Then, the master problem and feasibility check subproblems are solved iteratively.

*1) Bilinear Master Problem (**BMP**):* The master problem is to minimize the day-ahead operational cost with respect to the first-stage constraints (2)-(14) and Benders feasibility cuts (33). Given the binary indicator $z_n$ and the bilinear reformulation of chance constraint, we naturally develop Benders feasibility cuts in bilinear forms, as shown in (33). The bilinear master problem of $\tau$-th iteration is described as the following:

$$\min \sum_{t \in \mathcal{T}} \sum_{g \in \mathcal{G}} [C_g^{\text{NL}} u_{gt} + \sum_{s \in \mathcal{S}_g} C_{gs}^{\text{SU}} \delta_{gst} + C_g^{\text{SD}} y_{gt} + F_g(P_{gt})$$

$$+ C_g^{\text{R+}} R_{gt}^+ + C_g^{\text{R-}} R_{gt}^-] \quad (32)$$

$$\text{s.t.} \quad (2) - (14)$$

$$\text{Benders cuts shown at the bottom} \quad (33)$$

$$\sum_{n=1}^{N} \pi_n z_n \leq \epsilon, \ z_n \in \{0,1\}, \ \forall n \quad (34)$$

where $\hat{\Psi}_n^{(k)}$ indicates the optimal objective of feasibility check subproblem under the $n$-th scenario at the $k$-th iteration. $\hat{u}_{gt}^{(k)}$, $\hat{v}_{gt}^{(k)}$, $\hat{y}_{gt}^{(k)}$, $\hat{P}_{gt}^{(k)}$, $\hat{R}_{gt}^{+(k)}$ and $\hat{R}_{gt}^{-(k)}$ represent the fixed values of first-stage variables $u_{gt}$, $v_{gt}$, $y_{gt}$, $P_{gt}$, $R_{gt}^+$ and $R_{gt}^-$. Then, $\hat{\lambda}_{gtn}^{(k)}$, $\hat{\mu}_{gtn}^{(k)}$, $\hat{\rho}_{gtn}^{(k)}$, $\hat{\varphi}_{gtn}^{(k)}$, $\hat{\phi}_{gtn}^{(k)}$ and $\hat{\psi}_{gtn}^{(k)}$ represent the dual variables for constraints (45)-(50) of feasibility check subproblem.

Again, the bilinear Benders feasibility cuts (33) can be converted into linear forms using McCormick linearization method, see (35)-(41). Constraint (35) is shown at the bottom of this page. It can be observed that only the first stage decision variables $u_{gt}$, $v_{gt}$, $y_{gt}$, $P_{gt}$, $R_{gt}^+$ and $R_{gt}^-$ are involved in bilinear feasibility cuts (33). Then, the bilinear terms $u_{gt} z_n (= \tilde{u}_{gtn})$, $v_{gt} z_n (= \tilde{v}_{gtn})$, $y_{gt} z_n (= \tilde{y}_{gtn})$, $P_{gt} z_n (= \tilde{P}_{gtn})$, $R_{gt}^+ z_n (= \tilde{R}_{gtn}^+)$ and $R_{gt}^- z_n (= \tilde{R}_{gtn}^-)$ are linearized by including the following constraints:

$$0 \leq \tilde{u}_{gtn} \leq z_n; \quad u_{gt} + z_n - 1 \leq \tilde{u}_{gtn} \leq u_{gt} \quad (36)$$

$$0 \leq \tilde{v}_{gtn} \leq z_n; \quad v_{gt} + z_n - 1 \leq \tilde{v}_{gtn} \leq v_{gt} \quad (37)$$

$$0 \leq \tilde{y}_{gtn} \leq z_n; \quad y_{gt} + z_n - 1 \leq \tilde{y}_{gtn} \leq y_{gt} \quad (38)$$

$$\left\{ \hat{\Psi}_n^{(k)} + \sum_{t \in \mathcal{T}} \sum_{g \in \mathcal{G}} [\hat{\lambda}_{gtn}^{(k)}(u_{gt} - \hat{u}_{gt}^{(k)}) + \hat{\mu}_{gtn}^{(k)}(v_{gt} - \hat{v}_{gt}^{(k)}) + \hat{\rho}_{gtn}^{(k)}(y_{gt} - \hat{y}_{gt}^{(k)}) + \hat{\varphi}_{gtn}^{(k)}(P_{gt} - \hat{P}_{gt}^{(k)}) + \hat{\phi}_{gtn}^{(k)}(R_{gt}^+ - \hat{R}_{gt}^{+(k)}) \right.$$

$$\left. + \hat{\psi}_{gtn}^{(k)}(R_{gt}^- - \hat{R}_{gt}^{-(k)})] \right\} (1 - z_n) \leq 0 \qquad \forall n, k = 1, \ldots \tau - 1 \quad (33)$$

$$\left\{ \hat{\Psi}_n^{(k)} + \sum_{t \in \mathcal{T}} \sum_{g \in \mathcal{G}} \hat{\lambda}_{gtn}^{(k)}\left(u_{gt} - \hat{u}_{gt}^{(k)}\right) + \hat{\mu}_{gtn}^{(k)}\left(v_{gt} - \hat{v}_{gt}^{(k)}\right) + \hat{\rho}_{gtn}^{(k)}\left(y_{gt} - \hat{y}_{gt}^{(k)}\right) + \hat{\varphi}_{gtn}^{(k)}\left(P_{gt} - \hat{P}_{gt}^{(k)}\right) + \hat{\phi}_{gtn}^{(k)}\left(R_{gt}^+ - \hat{R}_{gt}^{+(k)}\right) + \hat{\psi}_{gtn}^{(k)}\left(R_{gt}^- - \hat{R}_{gt}^{-(k)}\right) \right\} -$$

$$\left\{ \hat{\Psi}_n^{(k)} z_n + \sum_{t \in \mathcal{T}} \sum_{g \in \mathcal{G}} \hat{\lambda}_{gtn}^{(k)}\left(\tilde{u}_{gtn} - \hat{u}_{gt}^{(k)} z_n\right) + \hat{\mu}_{gtn}^{(k)}\left(\tilde{v}_{gtn} - \hat{v}_{gt}^{(k)} z_n\right) + \hat{\rho}_{gtn}^{(k)}\left(\tilde{y}_{gtn} - \hat{y}_{gt}^{(k)} z_n\right) + \hat{\varphi}_{gtn}^{(k)}\left(\tilde{P}_{gtn} - \hat{P}_{gt}^{(k)} z_n\right) + \hat{\phi}_{gtn}^{(k)}\left(\tilde{R}_{gtn}^+ - \hat{R}_{gt}^{+(k)} z_n\right) \right.$$

$$\left. + \hat{\psi}_{gtn}^{(k)}\left(\tilde{R}_{gtn}^- - \hat{R}_{gt}^{-(k)} z_n\right) \right\} \leq 0; \qquad \forall n, k = 1, \ldots \tau - 1 \quad (35)$$

$$0 \leq \tilde{P}_{gtn} \leq P_g^{\max} z_n; \quad P_{gt} - P_g^{\max}(1-z_n) \leq \tilde{P}_{gtn} \leq P_{gt} \quad (39)$$
$$0 \leq \tilde{R}_{gtn}^+ \leq R_g^{+\max} z_n; \quad R_{gt}^+ - R_g^{+\max}(1-z_n) \leq \tilde{R}_{gtn}^+ \leq R_{gt}^+ \quad (40)$$
$$0 \leq \tilde{R}_{gtn}^- \leq R_g^{-\max} z_n; \quad R_{gt}^- - R_g^{-\max}(1-z_n) \leq \tilde{R}_{gtn}^- \leq R_{gt}^- \quad (41)$$

**Remark:** Again, due to the fact that the upper bounds of binary variables $u_{gt}, v_{gt}, y_{gt}$ and continuous variables $P_{gt}, R_{gt}^+, R_{gt}^-$ are naturally available, McCormick linearized constraints in (35)-(41) define the feasible set, whose projection is just the one defined by the bilinear constraint (33). They are not only computationally equivalent but also theoretically equivalent.

Note that the linearized bilinear master problem becomes an MIP and it can be solved by commercial solvers. Through solving the master problem, the first-stage decision variables are obtained, i.e., $\hat{u}_{gt}^{(\tau)}, \hat{v}_{gt}^{(\tau)}, \hat{y}_{gt}^{(\tau)}, \hat{P}_{gt}^{(\tau)}, \hat{R}_{gt}^{+(\tau)}$ and $\hat{R}_{gt}^{-(\tau)}$. Then, they are passed on to all feasibility check subproblems, see (45)-(50).

*2) Feasibility Check Subproblem for Each Scenario* (**FSP$_n$**): This subproblem is used to check whether the current commitment of thermal units can accommodate the fluctuation of wind power generation and load consumption in each scenario. The feasibility check subproblem of the $n$-th scenario at the $\tau$-th iteration is described as the following:

$$\min \ \Psi_n^{(\tau)} = \sum_{t \in \mathcal{T}} \left( \eta_{tn}^{1+} + \eta_{tn}^{1-} + \sum_{l \in \mathcal{L}} \eta_{ltn}^2 \right) \quad (42)$$

s.t. $(17)-(23)$

$$\sum_{g \in \mathcal{G}} P_{gtn} + \sum_{q \in Q} W_{qtn} = \sum_{b \in \mathcal{B}} L_{btn} + \eta_{tn}^{1+} - \eta_{tn}^{1-} \quad \forall t \quad (43)$$

$$-C_l - \eta_{ltn}^2 \leq \sum_{b \in \mathcal{B}} K_{lb} \left( \sum_{g \in \mathbf{G}(b)} P_{gtn} + \sum_{q \in \mathbf{W}(b)} W_{qtn} - L_{btn} \right)$$
$$\leq C_l + \eta_{ltn}^2 \quad \forall l, t \quad (44)$$

$$u_{gt} = \hat{u}_{gt}^{(\tau)} \ : \ \lambda_{gtn}^{(\tau)} \quad \forall g, t \quad (45)$$
$$v_{gt} = \hat{v}_{gt}^{(\tau)} \ : \ \mu_{gtn}^{(\tau)} \quad \forall g, t \quad (46)$$
$$y_{gt} = \hat{y}_{gt}^{(\tau)} \ : \ \rho_{gtn}^{(\tau)} \quad \forall g, t \quad (47)$$
$$P_{gt} = \hat{P}_{gt}^{(\tau)} \ : \ \varphi_{gtn}^{(\tau)} \quad \forall g, t \quad (48)$$
$$R_{gt}^+ = \hat{R}_{gt}^{+(\tau)} \ : \ \phi_{gtn}^{(\tau)} \quad \forall g, t \quad (49)$$
$$R_{gt}^- = \hat{R}_{gt}^{-(\tau)} \ : \ \psi_{gtn}^{(\tau)} \quad \forall g, t \quad (50)$$

In order to ensure that the second stage constraints (15)-(23) are feasible, we relax power balance constraint (15), transmission line capacity constraint (16) by introducing non-negative slack variables $\eta_{tn}^{1\pm}$ and $\eta_{ltn}^2$. The objective of **FSP$_n$** is to minimize the sum of all slack variables. If the objective $\Psi_n^{(\tau)}$ is larger than a preset threshold (it is set to 0 in this paper), which means that constraints (15)-(23) are infeasible under the $n$-th scenario, then the (linearized) bilinear feasibility cut (35) will be generated and added into the master problem **BMP**.

*D. Summary*

Using the bilinear Benders decomposition algorithm to solve the proposed chance-constrained UC model works as the following. The flow chart is presented in Fig. 1.

**Step 1. Initialization:**

Initialize the iteration counter, $\tau = 0$. Set the lower bound $LB = -\infty$ and the upper bound $UB = +\infty$. Set the tolerance $\sigma$.

**Step 2. Iteration:**

**Step 2.a**: Update the iteration $\tau \leftarrow \tau + 1$. Solve the (linearized) bilinear master problem **BMP**.
(i) If it is infeasible, then terminate the algorithm.
(ii) Otherwise, obtain the optimal objective $\hat{V}_L^{(\tau)}$, the optimal commitment solution $\hat{u}_{gt}^{(\tau)}, \hat{v}_{gt}^{(\tau)}, \hat{y}_{gt}^{(\tau)}, \hat{P}_{gt}^{(\tau)}, \hat{R}_{gt}^{+(\tau)}, \hat{R}_{gt}^{-(\tau)}$ and the scenario indicator $\hat{z}_n^{(\tau)}$.
(iii) Update the lower bound $LB = \hat{V}_L^{(\tau)}$.

**Step 2.b**: For all scenarios whose $\hat{z}_n^{(\tau)} = 0$,
(i) Solve the feasibility check subproblem **FSP$_n$** and obtain the optimal objective $\hat{\Psi}_n^{(\tau)}$.
(ii) If $\hat{\Psi}_n^{(\tau)}$ is larger than zero, then supply the feasibility cut (35) to the master problem **BMP**.

**Step 2.c**: Construct an MIP problem for deriving the upper bound $UB$. This problem is composed of the objective function (1) and constraints (2)-(14), (16)-(23) and (27)-(31) with fixed values for the first-stage decision variables $u_{gt}, v_{gt}, y_{gt}, \delta_{gst}, P_{gt}, R_{gt}^+$ and $R_{gt}^-$ (obtained from Step 2.a). Then, solve this problem and obtain its optimal objective $\hat{V}_U^{(\tau)}$. If this problem is infeasible, $\hat{V}_U^{(\tau)}$ is set to positive infinity. Finally, update $UB = \min\{UB, \hat{V}_U^{(\tau)}\}$.

**Step 3. Stopping Condition:**

If $(UB - LB)/LB \leq \sigma$, then terminate with the commitment solution associated with the latest upper bound $UB$. Otherwise, return to Step 2.

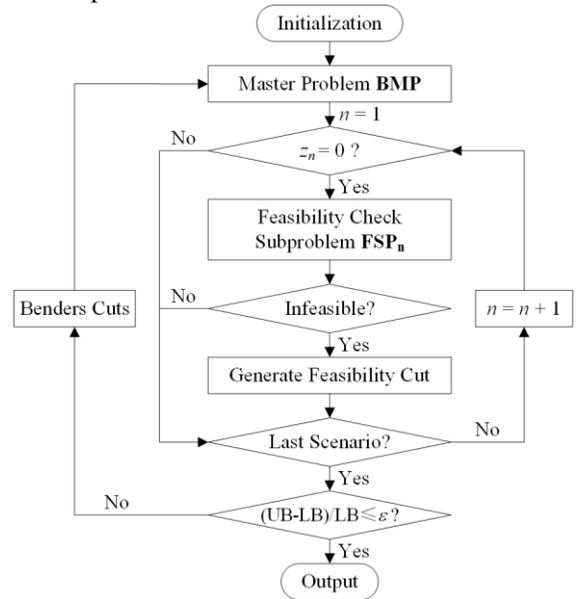

Fig. 1. Flow chart of the bilinear Benders decomposition algorithm.

## IV. COMPUTATIONAL RESULTS

A six-bus system and the modified IEEE 118-bus system are studied to validate the effectiveness of the proposed approach. The time horizon is 24 hours with one hour interval. Results of probabilistic wind power forecasts and probabilistic load forecasts come from GEFCom 2014 [48], and then they are used to generate several scenarios. The startup process of thermal unit is divided into two types, i.e., hot startup and cold startup [40]. Up/down-spinning reserve cost rates, $C_g^{R+}$ and $C_g^{R-}$, are set at 10% and 7%, respectively, of the linear coefficient of the quadratic fuel cost function [37]. The UC formulation and solution methodology are coded in GAMS environment using CPLEX 12.5. All experiments are implemented on a desktop computer with a processor clocking at 3.10 GHz and 8 GB of RAM.

### A. Six-Bus System

The six-bus system contains three thermal units, three loads and seven transmission lines. Detailed data of this testing system are given in [21, 36]. One wind farm is added to Bus 4. The percentage of wind power output with respect to the overall load for this system is 23.26%. This percentage is equal to the total forecasted wind power output divided by the total load demand over the 24-hour time horizon. In order to guarantee the fairness of comparison between Big-M and bilinear reformulations, we choose a small value for the Big-M parameter based on the summation of bounds of all variables involved in the constraint (24). As a result, Big-M is set to 1000 for the six-bus testing system. In this subsection, when using CPLEX to compute LP or MIP problems, the optimality gap $e$ is set to $10^{-4}$. For the bilinear Benders decomposition algorithm, the tolerance $\sigma$ is set to $10^{-4}$.

*1) Illustrative Example with Five Scenarios:* A rather small-scale case with only 5 scenarios is firstly studied to demonstrate the influence of chance constraints. The risk level $\epsilon$ is set to 20%. Table I shows the commitment solution of all thermal units obtained from chance-constrained UC models. According to the constraint (25), only one scenario at most is allowed to have load imbalance during the entire dispatching horizon. In other words, there is only one non-responsive scenario among all five scenarios. Experimental results show that the non-responsive scenario ($z_n = 1$) is Scenario 2.

On the other hand, it is easy to see that chance-constrained UC problems with fixed values of $z_n$ reduce to stochastic UC (SUC) problems. Then, we can compute SUC problems with different combinations of responsive and non-responsive scenarios. Table II lists the results of exhausting all possible SUC problems under the limitation of $\epsilon = 20\%$. From Table II, it can be seen that the operational cost reaches the lowest when Scenario 2 turns to be non-responsive, which is consistent with the results provided by chance-constrained UC models. Such results point out that chance-constrained UC models provide the solution with the minimal cost under the risk restriction on load imbalance.

TABLE I
UNIT COMMITMENT SOLUTION OF SIX-BUS SYSTEM

| Unit | Hours(1-24) |
|---|---|
| G1 | 111111 111111 111111 111111 |
| G2 | 000000 001111 111111 111110 |
| G3 | 000000 000001 111111 100000 |

TABLE II
OPERATIONAL COST OF SUC PROBLEMS WITH EXHAUSTING ALL COMBINATIONS OF SCENARIO INDICATORS $z_n$

| Non-responsive Scenario ($z_n = 1$) | Responsive Scenario ($z_n = 0$) | Operational Cost [$] |
|---|---|---|
| 1 | 2 3 4 5 | 80208 |
| 2 | 1 3 4 5 | 79606 |
| 3 | 1 2 4 5 | 80031 |
| 4 | 1 2 3 5 | 80208 |
| 5 | 1 2 3 4 | 80184 |

*2) Validation of Bilinear Reformulation and Benders Decomposition:* Detailed experiments are run to verify the effectiveness of the bilinear reformulation of chance constraint. The computational tractability of bilinear Benders decomposition algorithm is also validated in this subsection. Results are given in Tables III, IV and V. The stochastic UC (SUC) and chance-constrained UC (CC_UC) are compared in Table III. Notice that CC_UC models reduce to SUC models when the risk level $\epsilon$ is set to 0 (means that $z_n = 0$ for all $n$), given that the SUC model, which is adopted from [36], has no recourse cost. In Table III, several approaches to compute SUC and CC_UC models are specified as follows:

1) **CPX**: It denotes using CPLEX to directly solve the MIP formulation of SUC problems. SUC problems are formulated by (1)-(23) where the risk level $\epsilon$ is set to 0.
2) **CPX_BigM**: It represents using CPLEX to directly compute the Big-M reformulation of CC_UC problems [21-23]. The Big-M reformulation of CC_UC problem is formulated by (1)-(14), (24)-(25) and (16)-(23).
3) **CPX_Bilinear**: It denotes using CPLEX to directly compute the (McCormick linearized) bilinear reformulation of CC_UC problems. The bilinear reformulation of CC_UC problem is formulated by (1)-(14), (27)-(31) and (16)-(23).
4) **BD**: It represents using the bilinear variant of Benders decomposition algorithm to solve SUC/CC_UC problems. The number of iterations is displayed within the brackets in Table III.

The results shown in Table III are analyzed as the following. First, from Table III, it can be seen that solving CC_UC problems by CPLEX (columns "CPX_BigM" and "CPX_Bilinear") generally requires an order of magnitude more time than solving their corresponding SUC problems (column "CPX"). These results confirm that CC_UC models are much more difficult to solve than SUC models on condition of the same number of scenarios. Even for this small-scale six-bus testing system, a single chance constraint could lead to a significant increase of computational time for CC_UC problems. Second, it is found that the popular Big-M reformulation is not computationally efficient. Results in Table III indicate that solving the proposed bilinear reformulation of chance constraint (column "CPX_Bilinear") is faster than solving the Big-M reformulation (column "CPX_BigM"), verifying that the bilinear reformulation is stronger in computation. Third, the effectiveness of Benders decomposition is different for SUC and CC_UC problems. In SUC problems, directly solving UC problems by

CPLEX is typically more efficient than Benders decomposition. In CC_UC problems, however, Benders decomposition in the presented bilinear form is significantly faster than directly computing UC problems by CPLEX. From this result, we believe that the bilinear Benders decomposition algorithm is probably more appropriate to compute CC_UC problems than to compute SUC problems. Actually, it can be observed that the direct computation strategy using CPLEX is severely affected by chance constraints. In contrast, the presented bilinear variant of BD algorithm performs an order of magnitude faster. This observation is understandable because the total number of BD iterations does not increase significantly when chance constraints are imposed on SUC problems. Numerical results of BD algorithm show that CC_UC problems generally require just one or two more iterations to converge, in comparison with the corresponding SUC problems. Hence, the bilinear BD algorithm has strong capability to rapidly identify responsive scenarios from the scenario pool and simultaneously derive the solution of UC problems.

TABLE III
COMPUTATIONAL TIME FOR SIX-BUS SYSTEM WITH VARIED SCENARIO SIZES (RISK LEVEL OF CC_UC MODEL: 5%)

| Scenario Number | SUC | | CC_UC | | |
|---|---|---|---|---|---|
| | CPX [sec.] | BD [sec.] | CPX_BigM [sec.] | CPX_Bilinear [sec.] | BD [sec.] |
| 250 | 16.45 | 24.57(3) | 97.30 | 76.03 | 34.42(4) |
| 500 | 27.41 | 30.19(3) | 384.19 | 240.03 | 53.20(5) |
| 750 | 50.27 | 56.38(3) | 613.95 | 462.34 | 93.65(4) |
| 1000 | 94.27 | 140.33(4) | 1414.67 | 1041.34 | 202.42(5) |

Table IV shows numerical comparisons of integrality gap between the proposed bilinear reformulation and the conventional Big-M reformulation. In Table IV, column "MIP" demonstrates the objective value of MIP formulation (i.e., Big-M or bilinear reformulation) of CC_UC problems, and columns "CPX_BigM_Relax" and "CPX_Bilinear_Relax" give the objective value and integrality gap of the linear programming relaxations of Big-M and bilinear reformulations, respectively. Integrality gap (IG) is the ratio of objective values of MIP formulation and its relaxation. Results of Table IV indicate that the integrality gap of bilinear reformulation is consistently smaller than that of Big-M reformulation for different scenario numbers, verifying that the proposed bilinear reformulation is not only theoretically but also computationally stronger than the popular Big-M reformulation.

TABLE IV
INTEGRALITY GAP OF BILINEAR AND BIG-M REFORMULATIONS FOR CC_UC PROBLEMS (RISK LEVEL OF CC_UC MODEL: 5%)

| Scenario Number | MIP | CPX_BigM_Relax | | CPX_Bilinear_Relax | |
|---|---|---|---|---|---|
| | Obj.[$] | Obj.[$] | IG | Obj.[$] | IG |
| 250 | 80094.60 | 75822.70 | 1.0563 | 77209.04 | 1.0374 |
| 500 | 80151.78 | 76132.68 | 1.0528 | 77462.77 | 1.0347 |
| 750 | 80148.37 | 75833.56 | 1.0569 | 77275.04 | 1.0372 |
| 1000 | 80161.33 | 76323.13 | 1.0503 | 77490.90 | 1.0345 |

Table V gives computational time of using CPLEX directly or using BD algorithm to compute CC_UC problems with different risk levels. Column "[itr.]" presents the total number of BD iterations. The scenario number is set to 500. From Table V, it can be observed that the bilinear BD algorithm takes one order of magnitude less time than CPX_BigM/CPX_Bilinear on solving CC_UC problems. For CPX_BigM/CPX_Bilinear reformulations, the computational time is actually sensitive to the risk level. Note that the CC_UC problem with 10% risk level is the most difficult to compute by either of them, while the bilinear BD is rather robust against risk levels.

TABLE V
COMPUTATIONAL TIME OF CC_UC PROBLEMS FOR SIX-BUS SYSTEM WITH DIFFERENT RISK LEVELS (SCENARIO NUMBER: 500)

| Risk Level | CPX_BigM [sec.] | CPX_Bilinear [sec.] | BD [sec.] | [itr.] |
|---|---|---|---|---|
| 10% | 321.92 | 244.56 | 60.30 | 5 |
| 20% | 318.41 | 229.22 | 66.86 | 5 |
| 30% | 320.81 | 185.09 | 51.34 | 4 |
| 40% | 257.77 | 156.95 | 60.88 | 4 |

*3) Sensitivity Analysis of Different Risk Levels and Different Scenario Numbers:* Results of CC_UC problems with different risk levels are reported in Table VI. In this part, the bilinear BD algorithm, which has been testified to be the most efficient solution approach in Table III and V, is applied to solve CC_UC problems. From Table VI, it can be observed that the operational cost (column "objective") increases as the risk level decreases from 20% to 5%. It is understandable because the number of responsive scenarios increases with the decrease of risk level. As a result, the unit commitment solution should be more flexible to accommodate the fluctuation of wind power generation and load consumption in responsive scenarios, which may lead to a higher operational cost. The most extreme case is that the risk level is set to 0% and CC_UC problems reduce to SUC problems. In such a case, the operational cost is $80834. When the risk level is set to 100%, CC_UC problems reduce to deterministic UC problems. In this case, the operational cost is $79164, which is the smallest for all risk levels.

TABLE VI
COMPUTATIONAL RESULTS OF CC_UC PROBLEMS FOR SIX-BUS SYSTEM WITH DIFFERENT RISK LEVELS (SCENARIO NUMBER: 500)

| Risk Level | Objective [$] | Iteration | Time[sec.] |
|---|---|---|---|
| 5% | 80152 | 5 | 53.20 |
| 10% | 80138 | 4 | 44.90 |
| 15% | 79602 | 6 | 72.39 |
| 20% | 79414 | 4 | 42.40 |

In addition, we test the proposed model on different numbers of scenarios. Experimental results with 5% risk level are shown in Fig. 2. From Fig. 2, it can be seen that the objective function oscillates at the beginning when the scenario number is small and then converges slowly. When the number of scenarios is larger than 300, we observe that changes of the objective function value become less significant.

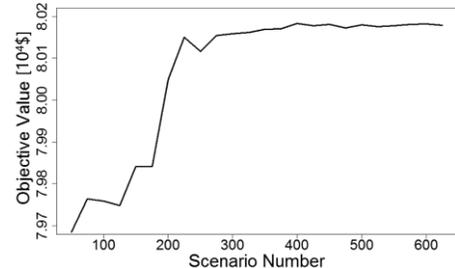

Fig. 2. Objective function of CC_UC problems with different scenario numbers.

## B. Modified 118-Bus System

The modified IEEE 118-bus system is applied to test the proposed approach for practical UC problems. This system has 54 thermal units, 186 transmission lines and 91 demands. The total peak load of 5516.08 MW occurs at hour 21. This system contains ten wind farms. The percentage of wind power output with respect to the total load for this system is 39.84%. Detailed data of 118-bus system are given in [49]. The Big-M parameter is set to 12,000. When using CPLEX to solve MIP problems, the optimality gap $e$ is set to 0.5%. To have a fair comparison, the tolerance of BD algorithm is also set to 0.5%. The limitation of computational time is set to 7,200 seconds.

Experimental results with different scenario numbers are reported in Table VII for comparison. The number of BD iterations is displayed within the brackets. In this subsection, because of the complexity of the 118-bus testing system, none of CPX, CPX_BigM or CPX_Bilinear is able to reach the optimality for any single instance within the preset time limit. Hence, the available objective and gap information before time limit are collected and presented as an alternative in Table VII. Note that CPX_Bilinear generally provides a better feasible solution with a smaller gap than that produced by CPX_BigM, which again confirms the strength of bilinear reformulation. Different from the aforementioned computational methods, our bilinear BD method computes the commitment solution of CC_UC problems with drastically less computational expenses. Moreover, comparing with the results in Table III, the number of Benders iterations does not have a large increase. Hence, we believe that the proposed BD method is less sensitive to the problem size and it has strong capability to solve practical CC_UC problems within reasonable time.

TABLE VII
COMPUTATIONAL RESULTS OF SUC AND CC_UC PROBLEMS FOR 118-BUS SYSTEM WITH DIFFERENT SCENARIO NUMBERS (RISK LEVEL: 5%)

| Scenario Number | SUC CPX | | CC_UC CPX_BigM | | CPX_Bilinear | | BD | |
|---|---|---|---|---|---|---|---|---|
| | Gap [%] | Obj. [$10^3$\$] | Gap [%] | Obj. [$10^3$\$] | Gap [%] | Obj. [$10^3$\$] | Obj. [$10^3$\$] | Time [sec.] |
| 100 | 1.52 | 1152.2 | 1.98 | 1162.9 | 1.59 | 1156.7 | 1138.0 | 721(5) |
| 200 | 4.07 | 1199.5 | 7.21 | 1246.1 | 3.96 | 1192.5 | 1131.2 | 768(5) |
| 300 | 1.53 | 1148.3 | 3.52 | 1172.6 | 1.32 | 1140.9 | 1120.4 | 575(6) |
| 400 | 1.86 | 1153.4 | 7.32 | 1235.8 | 1.91 | 1155.3 | 1126.8 | 743(7) |
| 500 | 2.39 | 1161.7 | 3.29 | 1181.0 | 2.50 | 1165.1 | 1128.6 | 1092(8) |

## V. CONCLUSION

This paper first presents a chance-constrained two-stage unit commitment (UC) model considering the uncertainty of load demand and wind power output. In this study, chance constraint guarantees the power balance being satisfied with a predefined high probability. Then, this UC model is reformulated as the bilinear mixed integer programming problem and then linearized by McCormick linearization method. Finally, to deal with the large number of random scenarios, a bilinear variant of Benders decomposition algorithm is developed to achieve fast computation for chance-constrained UC problems.

Numerical results on typical IEEE systems indicate that the bilinear formulation is computationally more efficient than the widely adopted Big-M based formulation. In particular, the bilinear Benders decomposition algorithm is an exact algorithm without any concern of Big-M issues and generally performs an order of magnitude faster than using a professional solver to directly compute both linear and bilinear chance-constrained UC models in terms of CPU time. In the future, more efforts would be made to develop appropriate enhancement techniques of bilinear Benders decomposition algorithm for further improving its computational efficiency.